\documentclass[12pt]{article}

\usepackage{amssymb}
\usepackage{pstricks}

\newcommand{\nc}{\newcommand}
\nc{\CC}{\mathbb C}
\nc{\NN}{\mathbb N}
\nc{\RR}{\mathbb R}
\nc{\ZZ}{\mathbb Z}

\renewcommand{\P}{\mathcal P}

\begin{document}

\title{Dirac eigenvalues and total scalar curvature}

\author{Bernd Ammann and Christian B\"ar}

\date{23. June 1999}
\maketitle

\begin{abstract}
\noindent
It has recently been conjectured that the eigenvalues $\lambda$ of the
Dirac operator on a closed Riemannian spin manifold $M$ of dimension $n\ge 3$
can be estimated from below by the total scalar curvature:
$$
\lambda^2 \ge \frac{n}{4(n-1)} \cdot \frac{\int_M S}{vol(M)}.
$$
We show by example that such an estimate is impossible.


\end{abstract}

\section{Introduction}
Let $M$ be a closed Riemannian spin manifold of dimension $n \ge 2$.
Then the eigenvalues of the Dirac operator $D$ acting on spinors form a
discrete subset of $\RR$ with finite multiplicities.
Let $S : M \to \RR$ denote the scalar curvature of $M$.
Let $\nabla$ be the Levi-Civita connection on spinors.
From the {\em Schr\"odinger-Lichnerowicz formula} \cite{schroedinger32a,
lichnerowicz63a}
\begin{equation}
D^2 = \nabla^\ast\nabla + \frac{S}{4}
\label{slf}
\end{equation}
it is clear that all eigenvalues $\lambda$ of $D$ must satisfy
\begin{equation}
\lambda^2 \ge \frac{1}{4}\min_{x\in M} S(x)
\label{schwach}
\end{equation}
because $\nabla^\ast\nabla$ is a nonnegative operator.
Of course, this is interesting only if the scalar curvature is positive
but then the estimate (\ref{schwach}) is never sharp.

The first sharp estimate was given by T.\ Friedrich who proved

{\bf Friedrich inequality \cite{friedrich80a}.}
{\em Let $M$ be a closed Riemannian spin manifold of dimension $n\ge 2$.
Then all eigenvalues $\lambda$ of the Dirac operator $D$ satisfy
$$
\lambda^2 \ge \frac{n}{4(n-1)}\min_{x\in M} S(x).
$$
}

This estimate is sharp in the sense that there are manifolds of positive
scalar curvature where equality is attained for the eigenvalue of smallest
modulus.
The round sphere is the simplest such example.
Friedrich's proof is based on a modification of the Levi-Civita
connection and a formula similar to (\ref{slf}).

Still, the Friedrich inequality is nontrivial only if the scalar curvature
is strictly positive.
To allow at least some negative scalar curvature it can be improved in 
two ways.
The first one works if the dimension is $n\ge 3$.

{\bf Hijazi inequality \cite{hijazi86a}.}
{\em Let $M$ be a closed Riemannian spin manifold of dimension $n\ge 3$.
Then all eigenvalues $\lambda$ of the Dirac operator $D$ satisfy
$$
\lambda^2 \ge \frac{n}{4(n-1)}\mu_1
$$
where $\mu_1$ is the first eigenvalue of the Yamabe operator
$Y = 4\cdot\frac{n-1}{n-2}\cdot\Delta + S$.
}

Since $\Delta$ is a nonnegative operator the Hijazi inequality implies the
Friedrich inequality if $n\ge 3$.
Hijazi's proof uses Friedrich's modification of the Levi-Civita connection
and a clever conformal change of the metric.
In dimension $n=2$ one can use a more complicated modification of the 
Levi-Civita connection and the Gauss-Bonnet theorem to show

{\bf B\"ar inequality \cite{baer91a,baer92d}.}
{\em Let $M$ be a closed and connected Riemannian spin manifold of dimension 
$n=2$.
Then all eigenvalues $\lambda$ of the Dirac operator $D$ satisfy
$$
\lambda^2 \ge \frac{2\pi\cdot\chi(M)}{area(M)}
$$
where $\chi(M)$ is the Euler number of $M$.
}

The estimate is nontrivial only if the surface is of genus 0 but
this does not come as a surprise because surfaces of genus $\ge 1$
are know to have the Dirac eigenvalue 0 for suitable choice of spin
structure and Riemannian metric.

It is now tempting to search for an improvement of the Hijazi inequality
and to conjecture the following

{\bf Conjecture \cite{friedrich-kim99a}.}
{\em Let $M$ be a closed Riemannian spin manifold of dimension $n\ge 2$.
Then all eigenvalues $\lambda$ of the Dirac operator $D$ satisfy
\begin{equation}
\lambda^2 \ge \frac{n}{4(n-1)}\frac{\int_M S}{vol(M)}.
\label{conj}
\end{equation}
}

In dimension 2 this is exactly the B\"ar inequality. 
Unfortunately, the conjecture is false for $n \ge 3$ as we shall see.
We will prove

{\bf Theorem.}
{\em Let $M$ be a closed spin manifold of dimension $n\ge 3$.
Then there exist constants $0 < C_1 \le C_2 \le \ldots$ and a set $\P$ of 
Riemannian metrics on $M$ such that
\begin{itemize}
\item
The $k^{th}$ Dirac eigenvalue $\lambda_k(g)$ (ordered by magnitude of its
modulus) for the Riemannian metric $g$ is bounded by $C_k$,
$$
\lambda_k(g)^2 \leq C_k
$$
for all $g \in \P$.
\item
The normalized total scalar curvature for the Riemannian metric $g$ is 
unbounded from above,
$$
\sup_{g\in\P}\frac{\int_M S_g}{vol_g(M)} = \infty .
$$
\end{itemize}
}

Therefore an estimate of the type
$$
\lambda_k^2 \ge C(k,n)\cdot\frac{\int_M S}{vol(M)}
$$
is impossible in dimension $n\ge 3$.

\section{Pinocchio Metrics}
In this section we will prove the theorem by explicitly constructing
the set of metrics $\P$ with the desired properties.
Let $M$ be a closed differentiable manifold of dimension $n\ge 3$ and
a fixed spin structure.
We choose a Riemannian metric $g_0$ on $M$ such that $(M,g)$ contains an
embedded Euclidean ball $B$ of radius 1.

\begin{center}
\psset{unit=1mm}
\psset{linewidth=1pt}
\begin{pspicture}(110,90)
\pscurve(35,19)(45,11.5)(50,8.5)(60,7)(68,8)(70,9)(73,10)(78,20)(79,25)(77,33)(74,34)(73,37)(73,40)(74,50)(74,58)(74,59)(75,61)(77,63)(80,65)(83,70)(80,78)(60,84)(40,80)(30,77)(20,73)(10.5,60)(10,50)(11.5,42)(12,40)(11.5,37)(10.5,30)(12,25)(20,20)(24,19.5)(28,20.5)

\psecurve(64,20)(65,16.5)(66,15)(69,14)(72,15)(73.5,18.5)(74,21)(74.3,25)
\pscurve(66,15)(67,17)(70,19)(73.5,18.5)
\pscurve(62,45)(62.8,41)(67,39)(70,40)(71,43)(71.2,45)(71.2,48)
\pscurve(62.7,42)(66,49)(69,48.5)(71.2,45)
\pscurve(51,39)(51.8,35)(56,33)(59,34)(60,37)(60.2,39)(60.2,42)
\pscurve(51.7,36)(55,43)(58,42.5)(60.2,39)

\pscircle(67,28){6}

\pscurve(38,22)(35,20)(32,19.2)(25.5,25)(27,30)(35,32)
\pscurve(35,23)(31.5,23.5)(30,25)(31,27)(33,28)
\pscurve(31.9,27.2)(33.5,25)(33,23.2)

\rput(33,60){$M$}

\rput(67,28){$B$} 

\end{pspicture}
\end{center}
\begin{center}
\bf Fig.\ 1
\end{center}

Write the Euclidean ball $B$ as a union of two annuli and one smaller 
ball, $B = A_1 \cup A_2 \cup A_3$, where $A_1 = \{ x \in \RR^n\ |\ 2/3 \leq 
|x| \leq 1 \}$, $A_2 = \{ x \in \RR^n\ |\ 1/3 \leq |x| \leq 2/3 \}$ and
$A_3 = \{ x \in \RR^n\ |\ |x| \leq 1/3 \}$.
Now fix two parameters $0 < r < 1$ and $L> 0$.
Choose a Riemannian metric $g_{r,L}$ on $M$ with the following properties:
\begin{itemize}
\item
$g_{r,L}$ coincides with $g_0$ on $M-B$
\item
$g_{r,L}$ is independent of $L$ on  $A_1$ and on $A_3$
\item
$(A_2,g_{r,L})$ is isometric to $S^{n-1}(r) \times [0,L]$ with the product
metric where $S^{n-1}(r)$ denotes the round sphere of constant sectional
curvature $1/r^2$.
\end{itemize}

\begin{center}
\psset{unit=1mm}
\psset{linewidth=1pt}
\begin{pspicture}(110,90)
\pscurve(35,19)(45,11.5)(50,8.5)(60,7)(68,8)(70,9)(73,10)(78,20)(79,25)(77,33)(74,34)(73,37)(73,40)(74,50)(74,58)(74,59)(75,61)(77,63)(80,65)(83,70)(80,78)(60,84)(40,80)(30,77)(20,73)(10.5,60)(10,50)(11.5,42)(12,40)(11.5,37)(10.5,30)(12,25)(20,20)(24,19.5)(28,20.5)

\psecurve(64,20)(65,16.5)(66,15)(69,14)(72,15)(73.5,18.5)(74,21)(74.3,25)
\pscurve(66,15)(67,17)(70,19)(73.5,18.5)
\pscurve(62,45)(62.8,41)(67,39)(70,40)(71,43)(71.2,45)(71.2,48)
\pscurve(62.7,42)(66,49)(69,48.5)(71.2,45)
\pscurve(51,39)(51.8,35)(56,33)(59,34)(60,37)(60.2,39)(60.2,42)
\pscurve(51.7,36)(55,43)(58,42.5)(60.2,39)

\pscircle(67,28){6}

\pscurve(38,22)(35,20)(32,19.2)(25.5,25)(27,30)(35,32)
\pscurve(35,23)(31.5,23.5)(30,25)(31,27)(33,28)
\pscurve(31.9,27.2)(33.5,25)(33,23.2)

\rput(33,60){$M$}


\rput{290}(67,28){
\scalebox{6}{
\psset{linewidth=.166pt}
\psframe*[linecolor=white](-.8,0)(.8,1)
\psecurve(-1.2,0)(-.9,0)(-.8,.3)(-.8,1)(-.8,1.5)
\psecurve(1.2,0)(.9,0)(.8,.3)(.8,1)(.8,1.5)
\psellipse(0,1)(.8,.4)
\psframe*[linecolor=white](-.8,1)(.8,6.3)
\psline(-.8,1)(-.8,6.3)
\psline(.8,1)(.8,6.3)
\psellipse(0,6.3)(.8,.4)
\psframe*[linecolor=white](-.8,6.3)(.8,6.7)
\psarc(0,6.3){.8}{0}{180}
}
}
\rput(67,28){$A_1$}
\rput(85,34.5){$A_2$}
\rput(103.4,41.6){$A_3$}

\end{pspicture}
\end{center}
\begin{center}
\bf Fig.\ 2
\end{center}

Heuristically, there is a nose of radius $r$ and length $L$ growing out of
the ball.
For this reason we call these metrics {\em Pinocchio metrics}\footnote{We
are indebted to A.\ Hornecker for suggesting this name.
There was another suggestion by L.\ Seeger to call them {\em Viagra metrics}
which we rejected with regard to our readers under age.}.
We set $\P := \{ g_{r,L}\ |\ 0 < r < 1,\ L> 0 \}$ and we check the 
properties required in the theorem.

{\bf Claim 1.}
{\em The $k^{th}$ eigenvalue $\lambda_k(r,L)^2$ of the square $D^2_{g_{r,L}}$ 
of the Dirac operator (w.r.t.\ the metric $g_{r,L}$) is bounded from above by 
a constant $C_k>0$ independent of $r$ and $L$.}

The proof is very simple.
Choose a $k$-dimensional vector space $V_k$ of spinors $\psi$ on $M$ vanishing
on $B$.
Then $\psi\in V_k$ can be considered a spinor for all metrics $g_{r,L}$.
We plug it into the Rayleigh quotient for the Dirac operator to get
\begin{eqnarray*}
\lambda_k(r,L)^2 &\leq& 
\sup_{\psi\in V_k, \psi\not= 0}
\frac{\int_M \langle D^2_{g_{r,L}}\psi,\psi\rangle_{g_{r,L}}}
{\int_M \langle \psi,\psi\rangle_{g_{r,L}}} \\
&=& 
\sup_{\psi\in V_k, \psi\not= 0}
\frac{\int_M \langle D^2_{g_0}\psi,\psi\rangle_{g_0}}
{\int_M \langle \psi,\psi\rangle_{g_0}} \\
&=:& C_k .
\end{eqnarray*}

{\bf Claim 2.}
{\em The normalized total scalar curvature is unbounded from above for $r\in 
(0,1)$ and $L \in (0,\infty)$.}

Let $\omega_k$ denote the volume of the $k$-dimensional unit sphere.
We compute
\begin{eqnarray*}
\frac{\int_M S_{g_{r,L}}}{vol_{g_{r,L}}(M)}
&=&
\frac{\int_{(M-B)\cup A_1\cup A_3} S_{g_{r,L}} + \int_{A_2} S_{g_{r,L}}}
{vol_{g_{r,L}}((M-B)\cup A_1\cup A_3)+vol_{g_{r,L}}(A_2)}\\
&=&
\frac{\int_{(M-B)\cup A_1\cup A_3} S_{g_{r,L}} 
+ L\cdot \frac{(n-1)(n-2)}{r^2}\cdot r^{n-1}\cdot\omega_{n-1}}
{vol_{g_{r,L}}((M-B)\cup A_1\cup A_3)
+ L\cdot r^{n-1}\cdot\omega_{n-1}}\\
&\longrightarrow&
\frac{(n-1)(n-2)}{r^2}
\end{eqnarray*}
for $L \to \infty$ because on $M-A_2$ the metric $g_{r,L}$ does not depend on 
$L$ by construction.
Hence 
$$
\sup_{L>0}\frac{\int_M S_{g_{r,L}}}{vol_{g_{r,L}}(M)} \ge 
\frac{(n-1)(n-2)}{r^2}
$$
and therefore
$$
\sup_{L>0,\ 0<r<1}
\frac{\int_M S_{g_{r,L}}}{vol_{g_{r,L}}(M)} = \infty .
$$
This proves the theorem.

\section{Concluding Remarks}

There is another reason why conjecture (\ref{conj}) cannot be true.
Assume that the closed spin manifold $M$ of dimension $n\ge 3$ has a metric
$g$ with harmonic spinors.
It is known that such a metric always exists if $n \equiv 0,1,3,7$ modulo 8
\cite{hitchin74a,baer96b}.
If the conjecture were true, then the total scalar curvature of every 
conformally equivalent metric $g_1$ would have to be nonpositive because
the multiplicity of the Dirac eigenvalue 0 is a conformal invariant.
But it is well-known that the total scalar curvature functional is not
bounded from above on any conformal class.
Here is the simple argument.
Write $g_1 = u^{4/(n-2)}\cdot g$ for some positive function $u$.
Then the scalar curvature of $g_1$ is given by 
$$
S_1 = u^{-\frac{n+2}{n-2}}\cdot Y(u)
$$
where $Y$ is the Yamabe operator as defined in the Hijazi inequality.
The volume element is 
$$
dvol_1 = u^{\frac{2n}{n-2}}dvol
$$
and thus the total scalar curvature is 
$$
\int_M S_1\ dvol_1 = \int_M u \cdot Y(u)\ dvol,
$$
cf.\ \cite[Ch.\ 1.J]{besse87a}.
Pick an eigenfunction $u$ of $Y$ for some positive eigenvalue $\mu > 0$.
The function $u$ will not be positive but for $\epsilon > 0$ we can
define $u_\epsilon := \sqrt{u^2 + \epsilon}$ and we see easily that
$$
\int_M u_\epsilon \cdot Y(u_\epsilon)\ dvol 
\stackrel{\epsilon \searrow 0}{\longrightarrow}
\int_M u \cdot Y(u)\ dvol 
= \mu \cdot \int_M u^2\ dvol
$$
Thus the total scalar curvature functional is not bounded from above on 
the conformal class of $g$.

For the proof of Claim 1 the use of the Dirac operator $D$ was rather 
inessential.
What we used is this:
The operator $D^2$ is elliptic and self-adjoint, it is computed locally out of the metric
and its derivatives, and its eigenvalues can be characterized
by variation of the Rayleigh quotient.
The theorem will still hold with $D^2$ replaced by any other operator
having these properties.
For example, we can take the Laplace operator acting on differential
forms.
If we denote the $k^{th}$ eigenvalue of the Laplace operator (w.r.t.\ the
Riemannian metric $g$) acting on $p$-forms by $\lambda^p_k(g)$, then we see
that 
$$
\sup_{g\in\P} \lambda^p_k(g) < \infty .
$$

Hence an estimate of the type
$$
\lambda^p_k \ge C(k,n)\cdot\frac{\int_M S}{vol(M)}
$$
is also impossible in dimension $n\ge 3$.



\providecommand{\bysame}{\leavevmode\hbox to3em{\hrulefill}\thinspace}

\vspace{1cm}

\parskip0ex

Mathematisches Institut

Universit\"at Freiburg

Eckerstr.~1

79104 Freiburg

Germany

\vspace{0.5cm}

E-Mail:
{\tt ammann@mathematik.uni-freiburg.de}

\hspace{1.4cm}
{\tt baer@mathematik.uni-freiburg.de}

WWW:
{\tt http://web.mathematik.uni-freiburg.de/home/ammann}

\hspace{1.4cm}
{\tt http://web.mathematik.uni-freiburg.de/home/baer}

\end{document}